\documentclass[a4paper,11pt]{amsart}

\usepackage[dvips]{graphicx}
\usepackage{amssymb}
\usepackage{amsmath}
\usepackage[latin1]{inputenc}

\newtheorem{theo}{Theorem}[section]

\newtheorem{lem}[theo]{Lemma}

\newtheorem{conj}[theo]{Conjecture}
\theoremstyle{definition}

\theoremstyle{remark} 

\theoremstyle{definition}

\newcommand{\cal}{\mathcal}

\begin{document}

\title[Diophantine properties of automatic real numbers]
{Diophantine properties of real numbers generated by finite automata}\author[B.~Adamczewski, J. Cassaigne]{Boris Adamczewski and Julien Cassaigne}

 %%%%%%%%%%%%%%%%%%%%%%%%%%%%%%%%%%%%%%%%%%%
\begin{abstract}
We study some diophantine properties of automatic real numbers and we present a method to derive irrationality measures for such numbers. 
As a consequence, we prove that the $b$-adic expansion of a Liouville number cannot be generated by a finite automaton, a conjecture due to Shallit.
\end{abstract}
%%%%%%%%%%%%%%%%%%%%%%%%%%%%%%%%%%%%%%%%%%%%%%%%%%%
\maketitle
%%%%%%%%%%%%%%%%
\section{Introduction}

The seminal work of Turing \cite{Turing} gives rise to a rough classification of real numbers. On one side we find computable real numbers, that is, real numbers whose binary (or more generally $b$-adic)  expansion can be produced by a Turing machine, while on the other side lie uncomputable real numbers which, in some sense, ``escape to computers". Though most real numbers belong to the second class (the first one being countable), classical mathematical constants are usually computable. By classical constants, we mean real numbers such as $\sqrt 2$ (or more generally algebraic real numbers), $\pi$ and $\zeta(3)$.  Note that the notion of period as considered by Kontsevitch and Zagier \cite{Kontsevitch_Zagier} could offer an interesting framework for (most of) these mathematical constants.  
Following the pioneering ideas of Turing, Hartmanis and Stearns \cite{Hartmanis_Stearns} 
proposed to emphasize the quantitative 
aspect of the notion of computability, and to take into
account the number $T(n)$ of operations  
needed by a (multitape) Turing machine to produce 
the first $n$ digits of the expansion. In this regard, a real number 
is considered all the more simple 
as its $b$-adic expansion can be produced very fast by 
a Turing machine. 
A general problem is then to determine where our  mathematical constants take place in such a classification.  It is a source of challenging open questions such as the Hartmanis--Stearns problem which asks whether there exists an irrational algebraic number computable in real-time, that is, with $T(n)=O(n)$. 

The present paper is partly motivated by the related question of how such a classification of real numbers based on Turing machines and theoretical computer science  fits into the one  
based on Diophantine approximation which was developped by Mahler \cite{Mahler32} during the 1930's. More modestly, we will focus on a special class of Turing machines of a particular interest: finite automata. They are one of the most basic models of computation and thus take place at the bottom of the hierarchy of Turing machines. In particular, such machines produce sequences in real-time. 
Roughly speaking, an infinite sequence ${\bf a}=(a_n)_{n\geq 0}$ is 
generated by a $k$-automaton if $a_n$ is a 
finite-state function of the base-$k$ 
representation of $n$. Consequently, automatic sequences share deep links with number theory as explained in detail in the book of Allouche and Shallit \cite{Allouche_Shallit}.

A real number is generated by a finite automaton (or in short automatic), if there exists a positive integer $b\geq 2$ such that its $b$-adic expansion can be ge\-nerated by a finite automaton. 
A particular attention was brought to these numbers since the following statement related to the Hartmanis--Stearns problem was suggested by Cobham \cite{Cobham68} in 1968: no irrational algebraic number can be generated by a finite automaton. Previously referred to as the Loxton--van der Poorten conjecture because of some attempts by these authors (see for instance \cite{LvdP1,LvdP2}), this result was recently proved in \cite{Adamczewski_Bugeaud, ABL}.
However it is still unknown whether constants such as $\pi$, $\zeta(3)$, or $e$ are or not automatic. 

Our main purpose is to  introduce a method to derive Diophantine pro\-perties for all automatic real numbers. This uses some of the ideas in \cite{Adamczewski_Bugeaud, ABL} together with classical techniques from Diophantine approximation and a careful combinatorial study of automatic sequences. In particular, we will prove that no Liouville number can be generated by a finite automaton; a result conjectured by Shallit \cite{Shallit}. Actually our approach is much more precise and it will provide an explicit general upper bound for the irrationality measure of any automatic real number.

\medskip

This article is organized as follows. We first present our main results, including the Shallit conjecture, in Section \ref{res}. In Section \ref{koksma}, we state a partial result towards a more general conjecture on automatic real numbers. This conjecture attributed to Becker involves Mahler's classification of real numbers and we will thus briefly recall some facts about it.  A  background on finite automata, words and morphisms can be found in Section \ref{auto}.  
Section \ref{preuves} is devoted to the proofs of our main results. As an application of our method, we derive an irrationality measure for the Thue--Morse--Mahler numbers in Section \ref{mahler}. Some comments and generalizations are gathered in Section \ref{remarks}. Finally, outlines of the proofs of these latter results end this paper.   
%%%%%%%%%%%%%%%%%%%%%%%%%%%%%%%%%%%%%%%%%%%%%%%%%%%%%%%%%%%%%%%%%%%%%%%%
\section{Main results}\label{res}

Diophantine approximation is essentially devoted to the estimate of  
the approximation of a given real number by rationals $p/q$, as a 
function of $q$.  A useful notion to tackle this question is 
the irrationality measure of a real number $\xi$, that we will 
denote by $\mu(\xi)$. It is defined as the supremum of the positive 
real numbers $\tau$ for which the inequality 
$$\left\vert \xi-\frac{p}{q}\right\vert<\frac{1}{q^{\tau}}$$
has infinitely many solutions $(p,q)\in{\mathbb Z}^2$, with relatively prime integers $p$ and $q$. 
Thus, $\mu(\xi)$ measures the quality of the best rational approximations to $\xi$. Note that 
it is in general a challenging problem to obtain an irrationality measure, 
i.e., to bound from above the irrationality measure of a given real number. 

Let us recall some well-known facts about this notion.  
First, the theory of continued fractions ensures that $\mu(\xi)\geq 2$, for 
any irrational number $\xi$; while $\mu(\xi)=1$ for a rational number $\xi$.  On the other hand, 
algebraic irrational numbers all have an irrationality measure equal to $2$, 
as follows from Roth's theorem \cite{Roth55}. Note that this is also the case for 
almost all real numbers (with respect to the Lebesgue measure), a 
result due to Khintchine \cite{Khintchine1} (see also \cite{Khintchine2}). 

At the opposite side, we find the Liouville numbers, introduced by Liouville in his famous 1844 paper \cite{Liouville}. Such numbers can be very well approximated by rationals. In more concrete terms,  
an irrational real number $\xi$ is called a Liouville numbers if it has an infinite irrationality measure, that is, if for any positive $\tau$ the inequality $$\left\vert \xi-\frac{p}{q}\right\vert<\frac{1}{q^{\tau}}$$
has at least one solution $(p,q)\in{\mathbb Z}^2$.

\medskip

Our main result  gives a positive answer to a conjecture of Shallit  
\cite{Shallit}.

\begin{theo}\label{liouville}
A Liouville number cannot be generated by a finite automaton.
\end{theo}

Actually, we will prove a quantitative version of Theorem \ref{liouville}. Indeed, our method provides the following explicit general upper bound for the irrationality measure of any number generated by a finite automaton. This bound only depends on three parameters which naturally appear in the study of automatic sequences and whose definitions are postponed to Section~\ref{auto}. 

\begin{theo}\label{mesir}
Let $k$ and $b$ be two integers at least equal to $2$ and let ${\bf a}=(a_n)_{n\geq 0}$ be an infinite sequence generated by a $k$-automaton and with values in $\{0,1,\ldots,b-1\}$. Let $m$ be the cardinality of the $k$-kernel of the sequence ${\bf a}$ and let $d$ be the cardinality of the internal  alphabet associated to ${\bf a}$.
Then, the irrationality measure $\mu(\xi)$ of the real number 
$$\xi:=\sum_{n=0}^{+\infty}\frac{a_n}{b^n}$$
satisfies
$$\mu(\xi)\leq dk(k^{m}+1).$$ 
\end{theo}

%%%%%%%%%%%%%%%%%%%%%%%%%%%%%%%%%%%%%%%%%%%%%%%%%%%%%%%%%%%%%%%%%%%%%%%%%
\section{The classification of Mahler}\label{koksma}

In 1932, Mahler \cite{Mahler32} introduced the first relevant classification of real numbers with respect to Diophantine approximation or more precisely with respect to their quality of approximation by algebraic numbers.  According to  Mahler's classification, real numbers are split into four classes, namely $A$-numbers, $S$-numbers, $T$-numbers, and $U$-numbers. We recall now how these four classes can be defined as first considered by Mahler. The reader is referred to \cite{Bugeaud} for a complete treatment of this topic.

\medskip

Let $n$ be a positive integer and $\xi$ be a real number. Then, we define $w_n(\xi)$ as the supremum of the real numbers $\omega$ for which there exist infinitely many integer polynomials $P(X)$ of degree at most $n$ and such that 
$$0<\vert P(\xi)\vert \leq H(P)^{-\omega},$$
where $H(P)$ denotes the height of the polynomial $P(X)$, that is, the maximum of the moduli of its coefficients.  
Note that according to this definition we have $w_1(\xi)=\mu(\xi)-1$. Then, we set 
$$w(\xi)=\limsup_{n\to\infty}\frac{w_n(\xi)}{n}\cdot$$
According to Mahler's classification, we say that $\xi$ is an 

\medskip

$A$-number, if $w(\xi)=0$;

\smallskip

$S$-number, if $0<w(\xi)<+\infty$;

\smallskip

$T$-number, if $w(\xi)=+\infty$ and $w_n(\xi)<+\infty$ for every $n\geq 1$;

\smallskip

$U$-number, if $w_n(\xi)=+\infty$ for some $n\geq 1$.

\medskip

\noindent  
In 1961, Wirsing \cite{Wirsing} proved that the class of $A$-numbers corresponds exactly to the set of real algebraic numbers.  It is also known since Mahler \cite{Mahler32b} that the set of $S$-numbers has full Lebesgue measure. Thus, most real numbers have to be $S$-numbers. However, it is generally rather delicate to prove that a given number lies in this class. 
The following conjecture was suggested by Becker in his correspondence with Shallit in 1993. It claims that irrational automatic real numbers should behave as do almost all real numbers with respect to Mahler's classification.

\begin{conj}\label{becker}
Automatic irrational real numbers are all $S$-numbers.
\end{conj} 

Since Liouville numbers are particular cases of $U$-numbers, this conjecture contains  
Theorem \ref{liouville}. A first step towards Conjecture \ref{becker} has recently been done in \cite{Adamczewski_Bugeaud} and \cite{ABL}. Their result can be reformulated as follows.

\begin{theo}[Adamczewski \& Bugeaud]\label{ABL}
The set of automatic irrational real numbers does not contains any $A$-number. 
\end{theo}

By combining ideas of the proof of Theorem \ref{mesir} with Theorem \ref{ABL} and a result of Baker \cite{Baker64}, we obtain the following result. 
 
\begin{theo}\label{unombres} 
Let $b\geq 2$ be an integer and ${\bf a}=(a_n)_{n\geq 0}$ be a non-eventually periodic automatic sequence. Let us assume that the  internal sequence associated with ${\bf a}$ begins in an overlap (see Section \ref{auto} for a definition). Then, the automatic irrational real number  
$$\xi:=\sum_{n=0}^{+\infty}\frac{a_n}{b^n}$$
is either a $S$-number or a $T$-number.
\end{theo}

We end this section by mentioning that very special examples of Conjecture \ref{becker} are already known. Indeed, for any positive integer $d\geq 2$, the real number 
$$\xi:= \sum_{n=0}^{+\infty}\frac{1}{2^{d^n}}$$
is an $S$-number, as follows from the work of Nishioka \cite{Nishioka}.

%%%%%%%%%%%%%%%%%%%%%%%%%%%%%%%%%%%%%%%%%%%%%%%%%%%%
\section{Finite automata, morphisms and Cobham's theorem}\label{auto}

Let $k \ge 2$ be an integer.
An infinite sequence ${\bf a}=(a_n)_{n\geq 0}$ is 
said to be $k$-automatic if $a_n$ is a finite-state function of the base-$k$ 
representation of $n$. This means that there exists a finite automaton starting 
with the $k$-ary expansion of $n$ as input and producing the term $a_n$ as 
output. A nice reference on this topic is the 
 book of Allouche and Shallit \cite{Allouche_Shallit}.

A more concrete definition of $k$-automatic sequences can be given as follows.
Denote by $\Sigma_k$ the set 
$\left\{0,1,\ldots,k-1\right\}$. 
By definition, a $k$-automaton is a $6$-tuple 
$$
A=\left(Q,\Sigma_k,\delta,q_0,\Delta,\tau\right),
$$
where $Q$ is a finite set of states, 
$\delta:Q\times\Sigma_k\rightarrow Q$ is the transition function, $q_0$ is 
the initial state, $\Delta$ is the output alphabet and $\tau : Q\rightarrow 
\Delta$ is the output function. For a state $q$ in $Q$ and for a finite 
word $W=w_1 w_2 \ldots w_n$ on the alphabet $\Sigma_k$, 
we define recursively $\delta(q,W)$ by 
$\delta(q,W)=\delta(\delta(q,w_1w_2\ldots w_{n-1}),w_n)$. 
Let $n\geq 0$ be an integer and let 
$w_r w_{r-1}\ldots w_1 w_0$ in $\left(\Sigma_k\right)^{r+1}$ 
be the $k$-ary expansion 
of $n$; thus, $n=\sum_{i=0}^r w_i k^{i}$. We denote by $W_n$ the 
word $w_0 w_1 \ldots w_r$. Then, a sequence ${\bf a}=(a_n)_{n\geq 0}$ is 
said to be $k$-automatic if there exists a $k$-automaton $A$ such that 
$a_n=\tau(\delta(q_0,W_n))$ for all $n\geq 0$.

A classical example of a $2$-automatic sequence is given by the 
Baum-Sweet sequence (see \cite{Baum_Sweet}) ${\bf a}=(a_n)_{n\geq 0}= 110110010100100110010 \ldots$
This sequence is defined 
as follows: $a_n$ is equal to $1$ if the 
binary expansion of $n$ contains no block of consecutive $0$'s of odd length, and $0$ otherwise. 
The Baum-Sweet sequence can be generated by the $2$-automaton
$$
A=\bigl(\{q_0, q_1,q_2\}, \{0, 1\}, \delta, q_0, \{0, 1\}, \tau \bigr),
$$
where 
$$\begin{array}{l}
\delta(q_0, 0) = q_1, \; \delta(q_0, 1) =q_0, \;  \delta (q_1, 0) =q_0, \;\delta(q_1,1)= q_2, \\ \\ \delta(q_2,0)=q_2, \;  \delta(q_2,1)=q_2, 
\end{array}$$
and $\tau (q_0) = \tau(q_1)=1$, $\tau (q_2) = 0$.

\bigskip

\noindent{\it Morphisms and Cobham's theorem}. 
For a finite set ${\cal A}$, we denote by ${\cal A}^*$ the free monoid 
generated by ${\cal A}$. The empty word $\varepsilon$ is the neutral element 
of ${\cal A}^*$. Let ${\cal A}$ and ${\cal B}$ be two finite sets. An 
application from ${\cal A}$ to ${\cal B}^*$ can be uniquely extended to a 
homomorphism between the free monoids ${\cal A}^*$ and ${\cal B}^*$. We 
call morphism from ${\cal A}^*$ to ${\cal B}^*$ such a 
homomorphism. If there is a positive integer
$k$ such that each element of ${\cal A}$ is mapped to a word of
length $k$, then the morphism is called $k$-uniform or simply
uniform. 
Similarly, an 
application from ${\cal A}$ to ${\cal B}$ can be uniquely extended to a  
homomorphism between the free monoids ${\cal A}^*$ and ${\cal
B}^*$. Such a homomorphism is called a coding (the term
`letter-to-letter' morphism is also used in the literature).

A morphism 
$\sigma$ from ${\cal A}^*$ into itself is said to 
be prolongable if there exists a 
letter $a$ such that $\sigma(a)=aW$, where the word $W$ is such 
that $\sigma^n(W)$ is a non-empty word for every $n\geq 0$. 
In that case, the 
sequence of finite words $(\sigma^n(a))_{n\geq 0}$ converges in 
${\cal A}^{\infty}={\cal A}^*\cup{\cal A}^{{\mathbb N}}$ endowed with its usual topology (see for instance \cite{Lothaire}, Chap. 2) 
to an infinite word ${\bf a}$ denoted $\sigma^{\infty}(a)$. This infinite word is clearly 
a fixed point for $\sigma$ (extended by continuity to infinite words) and we say that ${\bf a}$ 
is generated by the morphism 
$\sigma$. 

For instance, the Fibonacci morphism $\sigma$ defined over the alphabet 
$\{0,1\}$ by $\sigma(0)=01$ and $\sigma(1)=1$ is a non-uniform
morphism which generates the celebrated Fibonacci infinite 
word 
$$
{\bf a}=\lim_{n\to +\infty}\sigma^n(0)=010010100100101001\ldots
$$

Uniform morphisms and automatic sequences are strongly connected, as
shows the following result of Cobham \cite{Cobham72}. Theorem \ref{cob} in particular implies that finite 
automata produce sequences in real-time. 

\begin{theo}[Cobham]\label{cob}
A sequence is $k$-automatic if and only if it is
the image under a coding of a fixed point of a $k$-uniform morphism. 
\end{theo}

Thus, one  can always associate to a $k$-automatic sequence ${\bf a}$ a $5$-tuple $(\varphi,\sigma, i, {\cal A}, {\cal I})$, where $\sigma$ is a $k$-uniform morphism defined over a finite alphabet ${\cal I}$, $i$ is a letter of 
${\cal I}$, $\varphi$ is a coding from ${\cal I}$ into ${\cal A}$, and  such that 
$${\bf a}=\varphi({\bf i}),$$
with ${\bf i}=\sigma^{\infty}(i)$. The set ${\cal I}$ and the sequence ${\bf i}$ are respectively called the internal alphabet and the internal sequence associated to the $5$-uple $(\varphi,\sigma, i, {\cal A}, {\cal I})$. With a slight abuse of language, we will say in the sequel that ${\cal I}$ (resp. ${\bf i}$) is the internal alphabet (resp. internal sequence) associated to ${\bf a}$. Indeed, Cobham  gives in fact a canonical way to associate with ${\bf a}$ a $5$-uple $(\varphi,\sigma, i, {\cal A}, {\cal I})$.

\medskip

We recall now that the $k$-kernel of a sequence ${\bf a}=(a_n)_{n\geq 0}$ is defined as the set $N_k({\bf a})$ of all sequences $(a_{k^i\cdot n+j})_{n\geq 0}$, where $i\geq 0$ and $0\leq j<k^i$. This notion gives rise to 
another useful characterization of $k$-automatic sequences which was first proved in \cite{Eilenberg}.

\begin{theo}[Eilenberg]\label{eilenberg}
A sequence is $k$-automatic if and only if its $k$-kernel is finite.
\end{theo}

\bigskip

\noindent{\it Words and repetitive patterns}.  
We end this section with some notation about repetitions in Combinatorics on Words. Let ${\cal A}$ be a finite set.  
The length of a word
$W$ on the alphabet ${\cal A}$, that is, the number of letters
composing $W$, is denoted by $\vert W\vert$. For any positive integer $\ell$, we write
$W^\ell$ for the word $W\ldots W$ ($\ell$ times repeated concatenation
of the word $W$). More generally, for any positive real number
$x$, we denote by $W^x$ the word
$W^{\lfloor x\rfloor}W'$, where $W'$ is the prefix of
$W$ of length $\left\lceil(x- \lfloor x\rfloor)\vert W\vert\right\rceil$. 
Here, and in all what follows, $\lfloor y\rfloor$ and
$\lceil y\rceil$ denote, respectively, the integer part and the upper
integer part of the real number $y$. 
An overlap is a word of the form $WWa$ where $W$ is a non-empty finite word and $a$ is the first letter of $W$.
%%%%%%%%%%%%%%%%%%%%%%%%%%%%%%%%%%%%%%%%%%%%%%%%%%%
\section{Proof of the main results}\label{preuves}

This section is devoted to the proof of our main results, namely Theorems \ref{mesir} and \ref{unombres}.  
Before proving Theorem \ref{mesir}, an auxiliary result that we state just below is needed. 
Roughly speaking, this result will allow us   
 to control the repetitive patterns occurring as prefixes of any $k$-automatic sequence as a function of 
 the size of its $k$-kernel. This result will also be used in the proof of Theorem \ref{unombres}.

\begin{lem}\label{kp}
Let ${\bf u}$ be a non-eventually periodic $k$-automatic sequence  defined on an alphabet ${\cal A}$. Let $U\in{\cal A}^*$, $V\in{\cal A}^+$ and $s\in{\mathbb Q}$ be such that $UV^s$ is a prefix of the sequence ${\bf u}$. Let $m$ be the cardinality of the $k$-kernel of ${\bf u}$. Then, 
$$\frac{\vert UV^s\vert}{\vert UV\vert}<k^m.$$
\end{lem}

\begin{proof}
Let ${\bf u}=(u_n)_{n\geq 0}$ be a non-eventually periodic $k$-automatic sequence  defined on an alphabet ${\cal A}$. 
A triple $(h,p,l)$ of integers is said to be admissible (with respect to ${\bf u}$) if the following conditions hold:

\medskip

\begin{itemize}

\item[{\rm (i)}] $1\leq p\leq h\leq l$;

\smallskip\smallskip

\item[{\rm (ii)}] for all integers $n$ such that $h\leq n<l$, $u_{n-p}=u_n$;

\smallskip\smallskip

\item[{\rm (iii)}] $u_{l-p}\not=u_l$.

\end{itemize}

\medskip

Let $m$ be the cardinality of the set $N_k({\bf u})$. 
We shall prove that for every admissible triple $(h,p,l)$, $l<hk^m$. Let us  assume on the contrary that $(h,p,l)$ is an admissible triple such that $l\geq hk^m$ and we aim at deriving a contradiction. 

\medskip

For $i$ ranging from $0$ to $m$, let us construct a triple $(h_i,p,l_i)$ admissible with respect to a sequence ${\bf u}^{(i)}$ of $N_k({\bf u})$. 
We start with ${\bf u}^{(0)}={\bf u}$, $h_0=h$, and $l_0=l$. Then, given ${\bf u}^{(i)}=(u_n^{(i)})_{n\geq 0}$ and $(h_i,p,l_i)$, for $0\leq i<m$, let $r_i$ be the remainder in the division of $l_i$ by $k$. We define ${\bf u}^{(i+1)}=(u_{kn+r_i}^{(i)})_{n\geq 0}$ by extracting from ${\bf u}^{(i)}$ letters at positions congruent to $r_i$ modulo $k$. Furthermore, we set $l_{i+1}=\lfloor l_i/k\rfloor$ and $h_{i+1}=l_{i+1}+p-\lfloor (l_i+p-h_i)/k\rfloor$. 

\medskip

We are now going to prove by induction on $i$ the following:

\medskip

\begin{itemize}
\item[{\rm (a)}] ${\bf u}^{(i)}\in N_k({\bf u})$;

\smallskip\smallskip

\item[{\rm (b)}] $h_i\leq h_0$;

\smallskip\smallskip

\item[{\rm (c)}] $l_i\geq h_0k^{m-i}$;

\smallskip\smallskip

\item[{\rm (d)}] $(h_i,p,l_i)$ is admissible with respect to ${\bf u}^{(i)}$.

\end{itemize}

\medskip
 
 By assumption, this assertion is true for $i=0$. Let us assume that it holds for some $i$, $0\leq i<m$, and let us prove that this also is the case for $i+1$.  
  
$\rm (a)$  It is first clear that ${\bf u}^{(i)}\in N_k({\bf u})$ implies that ${\bf u}^{(i+1)}\in N_k({\bf u})$. 

$\rm (b)$ We now prove the useful fact that  $h_{i+1}\leq h_i$. Indeed, we have 
 $$\begin{array}{ll}
 h_i-h_{i+1}&=h_i-p-l_{i+1}+\left\lfloor (l_i+p-h_i)/k\right\rfloor\\ \\
 &=\left\lfloor \frac{k(h_i-p-l_{i+1})+l_i+p-h_i}{k}\right\rfloor\\ \\
 &=\left\lfloor \frac{r_i+(k-1)(h_i-p)}{k}\right\rfloor\geq 0\, ,
 \end{array}$$
 since $r_i\geq 0$, $k\geq 2$ and by assumption $h_i\geq p$. Consequently, $h_{i+1}\leq h_i\leq h_0$ and $\rm (b)$ holds.
 
$\rm (c)$ Since by assumption $l_i\geq h_0k^{m-i}$, we have that $l_{i+1}=\lfloor l_i/k\rfloor\geq \lfloor (k^{m-i}h_0)/k\rfloor=k^{m-i-1}h_0$.  

 $\rm (d)$  Now, let us prove that $(h_{i+1},p,l_{i+1})$ is admissible with respect to ${\bf u}^{(i+1)}$. As $p-h_i\leq 0$, $p-h_{i+1}=\lfloor (l_i+p-h_i)/k\rfloor-l_{i+1}\leq \lfloor l_i/k\rfloor-l_{i+1}=0$. Moreover, we get from $(c)$ that 
 $$
 l_{i+1}-h_{i+1}
 =\left\lfloor (l_i+p-h_i)/k\right\rfloor-p\geq\left\lfloor \frac{h_0k^{m-i}+p-h_i}{k}\right\rfloor-p\, ,$$
and since $h_i\leq h_0$ by $(b)$ and $i<m$, we obtain
 $$\begin{array}{ll}
 l_{i+1}-h_{i+1} &\geq\left\lfloor \frac{h_ik^{m-i}+p-h_i}{k}\right\rfloor-p
=\left\lfloor \frac{h_i(k^{m-i}-1)-(k-1)p}{k}\right\rfloor 
\\ \\
&\geq  \left\lfloor \frac{(h_i-p)(k-1)}{k}\right\rfloor\geq 0.
 \end{array}$$
 We thus have 
 \begin{equation}\label{i}
 1\leq p\leq h_{i+1}\leq l_{i+1}.
 \end{equation}
 
 By assumption, we have $u^{(i)}_{n-p}=u^{(i)}_n$, for every integer $n$ such that $h_i\leq n<l_i$, and thus, 
 \begin{equation}\label{ineg1}
 u^{(i)}_{n-jp}=u^{(i)}_n, \end{equation}
 for all integers $n$ and $j$ such that $j\geq 0$ and $h_i+(j-1)p\leq n<l_i$. 
  On the other hand, we get that $$kh_{i+1}+r_i=l_i+kp-k\lfloor (l_i+p-h_i)/k\rfloor\geq h_i+(k-1)p.$$ Thus, 
 $u^{(i)}_{n-kp}=u^{(i)}_n$, for  every integer $n$ such that $kh_{i+1}+r_i\leq n<l_i$,  
 and it {\it a fortiori} follows that $u^{(i)}_{k(n-p)+r_i}=u^{(i)}_{kn+r_i}$, for every integer $n$ such that  $kh_{i+1}+r_i\leq kn+r_i<kl_{i+1}+r_i$. This implies that
 \begin{equation}\label{ii}
 u^{(i+1)}_{n-p}=u^{(i+1)}_n, 
  \end{equation}
 for every integer $n$ such that $h_{i+1}\leq n<l_{i+1}$.

 Since $kl_{i+1}+r_i=l_i$, we have $u^{(i+1)}_{l_{i+1}}=u^{(i)}_{l_i}$ and $u^{(i+1)}_{l_{i+1}-p}=u^{(i)}_{l_i-kp}$.  Moreover,  $l_{i+1}\geq h_{i+1}$ implies that $\lfloor (l_i+p-h_i)/k\rfloor\geq p$ and thus $l_i\geq h_i+(k-1)p$. Then, we obtain  $l_i-p\geq h_i+(k-2)p$ and  we infer from (\ref{ineg1}) that 
 $$u^{(i+1)}_{l_{i+1}-p}=u^{(i)}_{l_i-kp}=u^{(i)}_{(l_i-p)-(k-1)p}=u^{(i)}_{l_i-p}\not=u^{(i)}_{l_i}=u^{(i+1)}_{l_{i+1}}.$$
 Together with (\ref{i}) and (\ref{ii}), this proves that the triple $(h_{i+1},p,l_{i+1})$ is admissible with respect to the sequence ${\bf u}^{(i+1)}$, hence $\rm (d)$.

\medskip

We are now ready to end our proof by deriving a contradiction. Since $N_{k}({\bf u})$ contains $m$ elements, at least two of the sequences ${\bf u}^{(i)}$, $0\leq i\leq m$, have to be the same. There thus exist two integers $i$ and $j$, $0\leq i<j\leq m$, such that $(h_i,p,l_i)$ and $(h_j,p,l_j)$ are two admissible triples for the same sequence. In particular, we have $u^{(i)}_n=u^{(j)}_n$ for any non-negative integer $n$. Moreover, $l_j\geq h_0k^{(m-j)}\geq h_0\geq h_i$, while $l_j=\left\lfloor \frac{l_i}{l_{j-1}}\right\rfloor<l_i$. Thus, $h_i\leq l_j<l_i$, and $(ii)$ gives that $u^{(i)}_{l_j-p}=u^{(i)}_{l_j}$. On the other hand, $(iii)$ implies that $u^{(i)}_{l_j-p}=u^{(j)}_{l_j-p}\not=u^{(j)}_{l_j}=u^{(i)}_{l_j}$, a contradiction. If the word $UV^s$ is a prefix of ${\bf u}$, as ${\bf u}$ is non-eventually periodic, there exists a maximal $s'\in{\mathbb Q}$, $s'\geq s$, such that $UV^{s'}$ is a prefix of ${\bf u}$, and $(\vert UV \vert, \vert V\vert,\vert UV^{s'}\vert)$ is an admissible triple. Therefore,  $\vert UV^{s'}\vert<\vert UV\vert k^m$ and thus $\vert UV^s\vert<\vert UV\vert k^m$ since $s'\geq s$. This concludes the proof.
\end{proof}

\medskip

Here and after in this paper, we will consider rational numbers defined thanks to their $b$-adic expansion. Actually, we will often use the following construction: given two finite words $U$ and $V$ defined over the alphabet $\{0,1,\ldots,b-1\}$, we define the rational number 
$$p/q:=0.UVVV\ldots V\ldots,$$
having preperiod $U$ and period $V$ in its $b$-adic expansion. It may happen that $V=b-1$ and in such a case the $b$-adic expansion of $p/q$ is usually defined by $0.a_1\ldots(a_l+1)$, where $U=a_1a_2\ldots a_l$. In the sequel of the paper, we allow inproper $b$-adic expansions of rational numbers, that is expansions ending with $(b-1)\ldots(b-1)\ldots$. So, when writting $p/q:=0.a_1a_2\ldots a_n\ldots$, with $a_n$ ultimately equal to $b-1$, this will naturally mean $p/q=\sum_{n=0}^{+\infty}a_n/b^n$.

\medskip

Before proving Theorem \ref{mesir}, we need to state the following result. 

\medskip

\begin{lem}\label{dist}
Let $b\geq 2$ be an integer, $U$ and $V$ be two finite words defined over the alphabet $\{0,1,\ldots,b-1\}$ with length respectively equal to $r$ and $s$, and $p/q$ be a rational number with eventually periodic $b$-adic expansion 
$$\frac{p}{q}:=0.a_1a_2\dots=0.UVV\ldots V\ldots$$  Let $\xi:=0.b_1b_2\ldots$ 
be a real number such that there exists a positive integer $j>r$ satisfying:

\begin{itemize}
\item[\rm (i)] $a_n=b_n$, for $1\leq n< j$;

\smallskip

\item[\rm (ii)] $a_j\not=b_j.$
\end{itemize}
Then,
$$\left\vert\xi- \frac{p}{q}\right\vert>\frac{1}{b^{j+s}}\cdot$$
\end{lem}

\medskip

\begin{proof}
We set $a_j=l$ and $b_j=m$. We first assume that $l>m$. Then, $l$ is a positive integer and we have  
$$\frac{p}{q}>0.a_1a_2\ldots a_{j-1}l\underbrace{00\ldots 0\ldots0}_{\mbox{$s-1$ times }}l,$$
whereas 
$$\xi\leq0.a_1a_2\ldots a_{j-1}(m+1)\leq 0.a_1a_2\ldots a_{j-1}l.$$ 
This yields 
 $$\frac{p}{q}-\xi>\frac{1}{b^{j+s}}\cdot$$

Now, let us assume that $m>l$. Then,  we have $l<b-1$ and 
$$\begin{array}{ll}
p/q&<0.a_1a_2\ldots a_{j-1}l\underbrace{(b-1)(b-1)\ldots (b-1)}_{\mbox{$s-1$ times }}(l+1)\\\\
&\leq0.a_1a_2\ldots a_{j-1}l\underbrace{(b-1)(b-1)\ldots (b-1)}_{\mbox{$s$ times }},
\end{array}$$
whereas 
$$\xi\geq0.a_1a_2\ldots a_{j-1}m \geq 0.a_1a_2\ldots a_{j-1}(l+1).$$ This gives 
$$\xi- \frac{p}{q}>\frac{1}{b^{j+s}},$$
concluding the proof.
\end{proof}

\medskip

We are now ready to prove our main result.

\medskip

\begin{proof}[Proof of Theorem \ref{mesir}] Let $k$ and $b$ be two integers at least equal to $2$ and let ${\bf a}=(a_n)_{n\geq 1}$ be an infinite sequence generated by a $k$-automaton and with values in $\{0,1,\ldots,b-1\}$. 
Let $\xi$ be the automatic real number defined by 
$$\xi:=\sum_{n=1}^{+\infty}\frac{a_n}{b^n}\cdot$$
Without loss of generality, we can assume that ${\bf a}$ is not eventually periodic. Indeed, otherwise $\xi$ would be a rational number and we would have $\mu(\xi)=1$. 

Following Cobham's theorem (Theorem \ref{cob}), there thus exist a $k$-uniform morphism $\sigma$ defined over a finite alphabet ${\cal I}$, a letter $i$ of 
${\cal I}$ and a coding 
$\varphi$ from ${\cal I}$ into $\{0,1,\ldots,b-1\}$ such that 
$${\bf a}=\varphi({\bf i}),$$
where ${\bf i}=\sigma^{\infty}(i)$. Let us also denote by $d$ the cardinality of the finite set ${\cal I}$ and  by $m$ the cardinality of the $k$-kernel of the sequence ${\bf a}$.

\medskip

Let $\delta$ be a positive number and let $(p,q)$ be a pair of positive integers. 
We shall prove that,  for $q$  large enough, one always has 
\begin{equation}\label{mes}
\left\vert \xi -\frac{p}{q}\right\vert\geq \frac{1}{q^{M+\delta}},
\end{equation}
with $M=dk(k^m+1)$.  Then, we will deduce the following upper bound for the irrationality measure of $\xi$: 
$$\mu(\xi)\leq M,$$
as claimed in Theorem \ref{mesir}.

\medskip

We are first going to introduce a sequence of rational numbers converging quite quickly to $\xi$ (see Inequality (\ref{nsba})) and 
whose denominators do not grow too fast (see Inequality (\ref{cptv})). These rational numbers are obtained thanks to Cobham's theorem and they will play a central role in this proof. 

\medskip
 
It follows from the pigeonhole principle that there is a letter $a$ that occurs at least twice in the prefix of length $d+1$ of the sequence ${\bf i}$. There thus exists a prefix of the sequence ${\bf i}$ of the form $UaVa$, where $U$ and $V$ are (possibly empty) finite words with $\vert U\vert+\vert V\vert\leq d-1$. 
Then, for any non-negative integer $n$, the finite word $\sigma^n(UaVa)$ is a prefix of the sequence ${\bf i}$. This implies that the finite word $\varphi(\sigma^n(UaVa))$ is a prefix of the sequence ${\bf a}$. 
Let $U_n=\varphi(\sigma^n(U))$ and $V_n=\varphi(\sigma^n(aV))$. We set $r_n=\vert U_n\vert$, $s_n=\vert V_n\vert$, and $t_n=\vert\varphi(\sigma^n(a))\vert$ . Since $\sigma$ is a $k$-uniform morphism and $\varphi$ is a coding, we get $r_n=k^n\vert U\vert$, $s_n=k^n\vert aV\vert$, and $t_n=k^n$. We also set $q_n=b^{r_n}(b^{s_n}-1)$. An easy computation shows  that there exists a positive integer $p_n$ such that the rational number $p_n/q_n$ has the following eventually periodic $b$-adic expansion:
$$\frac{p_n}{q_n}:=0.U_nV_nV_nV_n\ldots V_n\ldots$$
or in short 
$$\frac{p_n}{q_n}=
0.U_n\overline{V_n}.$$
Since the word $\varphi(\sigma^n(UaVa))$ is a common prefix of the sequence ${\bf a}$ and of $U_n\overline{V_n}$, it follows that 
\begin{equation}\label{j}
\left\vert\xi-\frac{p_n}{q_n}\right\vert<\frac{1}{b^{r_n+s_n+t_n}}\cdot
\end{equation}
Moreover, 
$$\frac{t_n}{r_n+s_n}=\frac{k^n}{k^n\vert UaV\vert}=\frac{1}{\vert UaV\vert}\geq \frac{1}{d}.$$ This implies
\begin{equation}\label{nsba}
\left\vert\xi-\frac{p_n}{q_n}\right\vert<\frac{1}{q_n^{1+1/d}}\cdot
\end{equation}
On the other hand, $r_{n+1}=kr_n$ and $s_{n+1}=ks_n$. This implies that  
\begin{equation}\label{cptv}
q_{n+1}=b^{kr_n}(b^{ks_n}-1)<\left(b^{r_n}(b^{s_n}-1)\right)^{k+\varepsilon}=q_n^{k+\varepsilon},
\end{equation}
for any positive real number $\varepsilon$ and any integer $n$ large enough.  Now, we fix a positive real number $\varepsilon$ that we choose small enough to ensure that 
\begin{equation}\label{epsilon}
d\varepsilon+d\varepsilon^2+\varepsilon dk^m+\varepsilon dk<\delta/2.
\end{equation} 
 Let $n_0$ be an integer such that (\ref{cptv}) 
holds for any integer $n\geq n_0$.

\medskip

We have now to remark that Lemma \ref{kp} naturally gives rise to an upper bound for the approximation of $\xi$ by $p_n/q_n$. By Lemma \ref{kp}, 
if $U_nV_n^s$ is a prefix of the sequence ${\bf a}$, then $\vert U_nV_n^s\vert<\vert U_nV_n\vert k^m=
(r_n+s_n)k^m$.  
We deduce that $\xi$ and $p_n/q_n$ cannot have the same first $\vert U_nV_n\vert k^m$ digits. Following (\ref{j}), the conditions of Lemma \ref{dist} are thus fulfilled with $r_n+s_n+t_n<j\leq (r_n+s_n)k^m$.  
Consequently, we get that 
$$\left\vert \xi-\frac{p_n}{q_n}\right\vert\geq \frac{1}{b^{(r_n+s_n) k^m+s_n}}>
\frac{1}{\left(b^{r_n}(b^{s_n}-1)\right)^{ k^m+1+\varepsilon}}
=\frac{1}{q_n^{k^m+1+\varepsilon}},$$
for $n$ large enough.  There thus exists a positive integer $n_1$, $n_1\geq n_0$, such that 
\begin{equation}\label{mercikp}
\left\vert \xi-\frac{p_n}{q_n}\right\vert>\frac{1}{q_n^{k^m+1+\varepsilon}},
\end{equation}
for any $n\geq n_1$.

\medskip

We are now ready for the last step of the proof.  Let us consider a rational number $p/q$ with $q$ large enough and let us prove that Inequality (\ref{mes}) holds. Let us assume that $q$ satisfies the following three inequalities: 
$2q\geq q_{n_1+1}^{1/d}$, $q\geq 2^{1+2M/\delta}$ and $q\geq 2^{1+dk+\delta/2}$ . Moreover, let us assume that 
\begin{equation}\label{1}
 \left\vert\xi-\frac{p}{q}\right\vert<\frac{1}{(2q)^{1+d(k+\varepsilon)}}\cdot
 \end{equation}

 First, note that if (\ref{1}) is not satisfied, then (\ref{mes}) holds. Indeed, let us assume that (\ref{1}) is not satisfied by the rational $p/q$. We have $dk^{m+1}>2$ and (\ref{epsilon}) implies that $d\varepsilon<\delta/2$. We thus derive that   $2+d(k+\varepsilon)<M+\delta$,  
 and since $q\geq 2^{1+dk+\delta/2}> 2^{1+dk+d\varepsilon}$, we have 
 $$(2q)^{1+d(k+\varepsilon)}<q^{2+d(k+\varepsilon)}<q^{M+\delta}.$$ This ensures that Inequality (\ref{mes}) is satisfied by the rational $p/q$.

Now, since by assumption $2q\geq q_{n_1+1}^{1/d}$, it follows from (\ref{cptv}) that there exists a unique positive integer $n_2> n_1\geq n_0$ such that 
\begin{equation}\label{3}
q_{n_2-1}\leq (2q)^{d}<q_{n_2}<q_{n_2-1}^{k+\varepsilon}\cdot
\end{equation} 
The second triangular inequality gives 
\begin{equation}\label{2}
\left\vert \xi-\frac{p_{n_2}}{q_{n_2}}\right\vert\geq \left\vert\; \left\vert \xi-\frac{p}{q}\right\vert-\left\vert\frac{p_{n_2}}{q_{n_2}}-\frac{p}{q}\right\vert\;\right\vert\cdot
\end{equation}
If $p/q$ and $p_{n_2}/q_{n_2}$ are two distinct rationals, we obviously get that   
\begin{equation*}
 \left\vert\frac{p_{n_2}}{q_{n_2}}-\frac{p}{q}\right\vert\geq\frac{1}{qq_{n_2}}\cdot
 \end{equation*}
 Moreover, $q_{n_2}\leq q_{n_2-1}^{k+\varepsilon}\leq (2q)^{d(k+\varepsilon)}$ and we thus derive from the previous inequality  that 
 \begin{equation*}
 \left\vert\frac{p_{n_2}}{q_{n_2}}-\frac{p}{q}\right\vert\geq\frac{2}{(2q)^{1+d(k+\varepsilon)}}> 2 \left\vert \xi-\frac{p}{q}\right\vert\cdot
 \end{equation*}
Then, (\ref{1}) and (\ref{2}) imply that, if $p/q\not=p_{n_2}/q_{n_2}$,  
$$\left\vert \xi-\frac{p_{n_2}}{q_{n_2}}\right\vert>\frac{1}{2qq_{n_2}}\cdot$$
Now, we infer from (\ref{3}) that $2q<q_{n_2}^{1/d}$. This gives 
$$\left\vert \xi-\frac{p_{n_2}}{q_{n_2}}\right\vert>\frac{1}{q_{n_2}^{1+1/d}}\;,$$
which is a contradiction with (\ref{nsba}). 
It thus follows that, under the assumption (\ref{1}), we necessarily have $p/q=p_{n_2}/q_{n_2}$.  
In that case, since $n_2\geq n_1$, we infer from (\ref{mercikp}) that 
$$\left \vert \xi-\frac{p}{q}\right\vert=\left \vert \xi-\frac{p_{n_2}}{q_{n_2}}\right\vert\geq 
\frac{1}{q_{n_2}^{k^m+1+\varepsilon}}$$
and since $(2q)^{d(k+\varepsilon)}\geq q_{n_2}$, we obtain
$$\left \vert \xi-\frac{p}{q}\right\vert\geq\frac{1}{(2q)^{d(k+\varepsilon)(k^m+1+\varepsilon)}}
\cdot$$
We infer from (\ref{epsilon}) and from the assumption that $q\geq 2^{1+2M/\delta}$ that
 $$(2q)^{d(k+\varepsilon)(k^m+1+\varepsilon)}\leq (2q)^{M+\delta/2}\leq q^{M+\delta/2}2^{M+\delta/2}\leq q^{M+\delta}$$ 
and it implies  
$$\left \vert \xi-\frac{p}{q}\right\vert\geq\frac{1}{q^{M+\delta}}\cdot$$
Inequality (\ref{mes}) thus holds for any $$q\geq  \max\left\{\frac{q_{n_1+1}^{1/d}}{2},2^{1+2M/\delta},2^{1+dk+\delta/2}\right\}\cdot$$ This concludes the proof.
\end{proof}

We end this Section by proving Theorem \ref{unombres}. The proof follows the same line as the one of Theorem \ref{mesir}. The key point is a general criterion to prove that a real number is not a $U$-number. This result is due to Baker \cite{Baker64} and we recall it now.

\begin{theo}[A. Baker]\label{baker}
Let $\xi$ be a real number, $\varepsilon$ be a positive number and $(p_n/q_n)_{n\geq 1}$ be a sequence of distinct rational numbers. Let us assume that the following conditions hold:
$$\left\vert \xi-\frac{p_n}{q_n}\right\vert<\frac{1}{q_n^{2+\varepsilon}}$$
and 
$$\limsup_{n\to\infty} \frac{\log(q_{n+1})}{\log(q_n)}<+\infty\cdot$$
Then, $\xi$ is not a $U$-number.
\end{theo}

\begin{proof}[Proof of Theorem \ref{unombres}]
We keep the notation of Theorem \ref{mesir}. Let us assume that the sequence ${\bf i}$ begins in an overlap $UUa$, where $U$ is a finite word and $a$ is the first letter of $U$.  Then, for any non-negative integer $n$, the finite word $\varphi(\sigma^n(UUa))$ is a prefix of the sequence ${\bf a}$. 
Let $s_n=\vert\varphi(\sigma^n(U))\vert$. We also set $q_n=b^{s_n}-1$. There thus exists a positive integer $p_n$ such that the rational number $p_n/q_n$ has the following periodic $b$-adic expansion:
$$\frac{p_n}{q_n}=
0.\overline{\varphi(\sigma^n(U))}.$$
Since the word $\varphi(\sigma^n(UUa))$ is a prefix of the sequence ${\bf a}$, it follows that 
\begin{equation}\label{big}
\left\vert\xi-\frac{p_n}{q_n}\right\vert<\frac{1}{b^{2s_n+t_n}},
\end{equation}
where $t_n=\vert \varphi(\sigma^n(a))\vert$. Moreover, since $\sigma$ is a $k$-uniform morphism, $s_n=k^n\vert U\vert$ and $t_n=k^n$, and we get that 
\begin{equation}\label{nsba3}
\left\vert\xi-\frac{p_n}{q_n}\right\vert<\frac{1}{q_n^{2+\varepsilon}},
\end{equation}
where $\varepsilon=1/\vert U\vert$.
On the other hand, the definition of $q_n$ implies that   
\begin{equation}\label{cptv3}
\limsup_{n\to\infty} \frac{\log(q_{n+r})}{\log(q_n)}=k^r<+\infty,
\end{equation}
for any positive integer $r$.
We now infer from Lemmata \ref{kp} and \ref{dist} that 
\begin{equation}\label{sm}
\left\vert \xi-\frac{p_n}{q_n}\right\vert>\frac{1}{b^{s_n(k^m+1)}}\cdot
\end{equation}
There thus exists a positive integer $r$ large enough to ensure that the sequence 
$(p_{rn}/q_{rn})_{n\geq 0}$ consists of distinct rational numbers. Indeed, let $n'>n\geq 0$ be two integers such that $p_n/q_n=p_{n'}/q_{n'}$. It follows from (\ref{big}) and (\ref{sm}) that 
$2s_{n'}+t_{n'}<s_n(k^m+1)$, which gives $2k^{n'}\vert U\vert+k^{n'}<k^n\vert U\vert(k^m+1)$. We thus have  
$$k^{n'-n}<\frac{\vert U\vert(k^m+1)}{2\vert U\vert+1}<k^m$$
and $n'<n+m$.

Then, (\ref{nsba3}) and (\ref{cptv3}) show that the conditions of Theorem \ref{baker} are fulfiled. This implies that $\xi$ is not a $U$-number. On the other hand, we infer from Theorem \ref{ABL} that $\xi$ is not a $A$-number. This ends the proof.
\end{proof}

%%%%%%%%%%%%%%%%%%%%%%%%%%%%%%%%%%%%%%%%%%%%%%%
\section{An emblematic example: the Thue--Morse--Mahler numbers}\label{mahler}

In Theorem \ref{mesir}, we have obtained a general upper bound for the irrationality measure 
of real numbers generated by finite automata. It turns out that for a specific automatic real number the method introduced in the proof of Theorem \ref{mesir} will likely give rise to a better estimate. In the present Section, we choose to illustrate this general idea by considering as a particular example the Thue--Morse--Mahler numbers.

For an integer $b$ at least equal to $2$, we define the $b$-adic Thue--Morse--Mahler number $\xi_b$ by 
$$\xi_b=\sum_{k=0}^{+\infty}\frac{a_k}{b^k},$$
where the sequence ${\bf a}=(a_n)_{n\geq 0}$ is the Thue--Morse sequence. This sequence is defined as follows: $a_n$ is equal to $0$ (resp. to $1$)
if the sum of the digits in the 
binary expansion of $n$ is even (resp. is odd). Thus, the $b$-adic expansion of $\xi_b$ begins with  
$0.110100110010 \ldots$. These numbers were first considered by Mahler who proved that they are  transcendental in \cite{Mahler29} (see also \cite{Dekking}). 

It is easy to 
check that the Thue--Morse sequence can be generated by the $2$-automaton
$$
A=\bigl(\{q_0, q_1\}, \{0, 1\}, \delta, q_0, \{0, 1\}, \tau \bigr),
$$
where 
$$
\delta(q_0, 0) = \delta (q_1, 1) = q_0, \qquad
\delta(q_0, 1) = \delta (q_1, 0) = q_1,
$$
and $\tau (q_0) = 0$, $\tau (q_1) = 1$. It is also well-known, as a consequence of Theorem~\ref{cob}, that this sequence is the fixed point beginning with $0$ of the following binary $2$-uniform morphism:
$$\begin{array}{lcl}
& \sigma &\\
0 &\longmapsto &01\\
1& \longmapsto & 10
\end{array}$$
Theorem \ref{mesir} thus implies that the irrationality measure of $\xi_b$ satisfies 
$\mu(\xi_b)\leq 20$ (independently of the integer $b$), since in that case $d=k=m=2$ (the fact that $m=2$ is for instance proved in \cite{Allouche_Shallit}). However, a better use of the method described in the previous Section  
leads to the following improvement.

\begin{theo}\label{tmm}
For any $b\geq 2$, we have 
$$\mu(\xi_b)\leq 5.$$
\end{theo}

\begin{proof} From now on, we fix a positive integer $b\geq2$ and we set $\xi:=\xi_b/b$. Thus, 
$$
\xi=0.a_0a_1a_2a_3\ldots
$$ 
It is obvious that $\xi$ and $\xi_b$ have the same irrationality measure. 
Let $\delta$ be a positive number. We want to prove that the inequality
\begin{equation}\label{cont}
\left\vert\xi-\frac{p}{q}\right\vert<\frac{1}{q^{5+\delta}}
\end{equation}
has only finitely many solutions $(p,q)\in{\mathbb Z}^2$.

As in the proof of Theorem \ref{mesir}, we are first going to introduce an infinite sequence of rationals converging to $\xi$. 
We can remark that $\sigma^3(0)=01101001$ and thus the sequence ${\bf a}$ begins in the word $011010$. More generally, it follows that ${\bf a}$ begins 
in $\sigma^n(011010)$ for any non-negative integer $n$. Since $0110$ is a prefix of $\sigma^n(0)$ for $n\geq 2$, we deduce that the word 
\begin{equation}\label{deb}
\sigma^n(011)\sigma^n(01)0110=\left(\sigma^n(011)^{1+2/3}\right)0110
\end{equation}
is a prefix of ${\bf a}$ for $n\geq 2$.  
Now, we set $q_n=(b^{3\cdot 2^n}-1)$. It follows from an easy computation that there exists a positive integer $p_n$ such that the rational number $p_n/q_n$ has the following periodic $b$-adic expansion:
$$\frac{p_n}{q_n}:=0.\sigma^n(011)\sigma^n(011)\ldots\sigma^n(011)\ldots$$
or in short 
$$\frac{p_n}{q_n}=
0.\overline{\sigma^n(011)}.$$
The $b$-adic expansion of $p_n/q_n$ begins with  
$$\sigma^n(011)\sigma^n(011)=\sigma^n(011)\sigma^n(01)\sigma^n(1)$$
and we deduce that it begins with  
\begin{equation}\label{deb2}
\sigma^n(011)\sigma^n(01)1001,
\end{equation}
when $n\geq 2$. 
  Since $\sigma$ is a $2$-uniform morphism, we easily check that $\vert\sigma^n(W)\vert=2^n\vert W\vert$ for any finite word $W$. 
We thus infer from (\ref{deb}) and (\ref{deb2}) that 
the first $(5\cdot 2^n+1)$-th digits in the $b$-adic expansion of $p_n/q_n$ and of $\xi$ are the same while the following $4$ digits are respectively $1001$ and $0110$. 
 This implies 
\begin{equation*}
\frac{1}{b^{5\cdot 2^n+3}}\leq \left\vert \xi-\frac{p_n}{q_n}\right\vert< \frac{1}{b^{5\cdot 2^n+2}}\cdot
\end{equation*}  
Let $\varepsilon$ be a positive number chosen small enough to ensure that  
\begin{equation}\label{choix}
(2q)^{4+3\varepsilon/2}<(2q)^{5+11\varepsilon/2+3\varepsilon^2/2}<q^{5+\delta}
\end{equation} 
for any integer $q$ large enough, say for any $q\geq M>2$. The definition of $q_n$ implies that there 
exists a positive integer $n_0$ such that 
\begin{equation}\label{ap}
\frac{1}{q_n^{1+2/3+\varepsilon}}\leq \left\vert \xi-\frac{p_n}{q_n}\right\vert< \frac{1}{q_n^{1+2/3}}
\end{equation} 
for any $n\geq n_0$. 
On the other hand, we also have 
\begin{equation}\label{cptv2}
q_{n+1}=(b^{3\cdot2^{n+1}}-1)<q_n^{2+\varepsilon}
\end{equation}
for any integer $n$ large enough, say $n\geq n_1\geq n_0$.  

\medskip

Let us assume that there exists a rational number $p/q$ satisfying Inequality (\ref{cont}) 
and such that the following two inequalities hold: 
$2q\geq q_{n_1}^{2/3}$ and $q\geq M$.  We now aim at deriving a contradiction.

\medskip

Since by assumption $2q\geq q_{n_1}^{2/3}$, it follows from (\ref{cptv2}) that there exists a unique positive integer $n_2>n_1$ such that 
\begin{equation}\label{32}
q_{n_2-1}\leq (2q)^{3/2}<q_{n_2}<q_{n_2-1}^{2+\varepsilon}\leq 2q^{\frac{3}{2}\cdot(2+\varepsilon)}\cdot
\end{equation} 
The second triangular inequality gives 
\begin{equation}\label{22}
\left\vert \xi-\frac{p_{n_2}}{q_{n_2}}\right\vert\geq \left\vert\; \left\vert \xi-\frac{p}{q}\right\vert-\left\vert\frac{p_{n_2}}{q_{n_2}}-\frac{p}{q}\right\vert\;\right\vert\cdot
\end{equation}
Then, we infer from (\ref{cont}), (\ref{ap}) and (\ref{32}) that  $p/q$ and $p_{n_2}/q_{n_2}$ are two distinct rationals. Otherwise we would have from (\ref{ap}) that  
$$\left\vert\xi-\frac{p}{q}\right\vert=\left\vert\xi-\frac{p_{n_2}}{q_{n_2}}\right\vert\geq  \frac{1}{q_{n_2}^{1+2/3+\varepsilon}}$$
and then (\ref{32}) would give 
$$\left\vert\xi-\frac{p}{q}\right\vert\geq \frac{1}{(2q)^{3/2(2+\varepsilon)(1+2/3+\varepsilon)}}
= \frac{1}{(2q)^{5+11\varepsilon/2+3\varepsilon^2/2}}\cdot$$
In such a case, (\ref{choix}) would give  
$$\left\vert\xi-\frac{p}{q}\right\vert\geq\frac{1}{q^{5+\delta}}$$
and we would reach a contradiction with (\ref{cont}). 
Thus $p/q$ and $p_{n_2}/q_{n_2}$ are distinct and we immediatley get that   
\begin{equation*}
 \left\vert\frac{p_{n_2}}{q_{n_2}}-\frac{p}{q}\right\vert\geq\frac{1}{qq_{n_2}}\cdot
 \end{equation*}
We then derive, using (\ref{32}), (\ref{choix}), and (\ref{cont}), that
 \begin{equation*}
 \left\vert\frac{p_{n_2}}{q_{n_2}}-\frac{p}{q}\right\vert\geq\frac{2}{(2q)^{1+3/2(2+\varepsilon)}}
 =\frac{2}{(2q)^{4+3\varepsilon/2}}> \frac{2}{q^{5+\delta}}>2 \left\vert \xi-\frac{p}{q}\right\vert\cdot
 \end{equation*}
Then, (\ref{cont}) and (\ref{22}) imply that 
$$\left\vert \xi-\frac{p_{n_2}}{q_{n_2}}\right\vert>\frac{1}{2qq_{n_2}}\cdot$$
Now, we infer from (\ref{32}) that $2q<q_{n_2}^{2/3}$. This gives 
$$\left\vert \xi-\frac{p_{n_2}}{q_{n_2}}\right\vert>\frac{1}{q_{n_2}^{1+2/3}}\;,$$
which is a contradiction with (\ref{ap}). Thus, any rational $p/q$ satisfying  Inequa\-lity (\ref{cont})  
has  a denominator $q$ at most equal to $\max\{M,q_{n_1}^{2/3}/2\}$. 
This ends the proof. 

\end{proof}
%%%%%%%%%%%%%%%%%%%%%%%%%%%%%%%%%%%%%%%%%%%%%%
\section{Some generalizations and concluding remarks}\label{remarks}

A first possibility of extension of the present approach is to focus on the representations of real numbers in algebraic bases, that is, to replace the integer $b$ with an algebraic real number $\beta$ greater than $1$. We can consider both the analog of the Loxton--van der Poorten conjecture and of the Becker conjecture, the number field $\mathbb Q(\beta)$ taking the place of the field of rationals. For an algebraic real number $\beta$ and an automatic sequence $(a_n)_{n\geq 0}$ with values in $\{0,1,\ldots,\lfloor \beta\rfloor\}$, it is thus likely that the real number $$\xi:=\sum_{n\geq 0}\frac{a_n}{\beta^n}$$ either lies in $\mathbb Q(\beta)$ or is transcendental (this is the analog of the Loxton--van der Poorten conjectue), and is always an $S$-number in the latter case (this is the analog of Becker's conjecture). 

Actually, if we restrict our attention to the $\beta$-expansions introduced by R\'enyi \cite{Renyi} and on a Pisot or a Salem base, the analog of the Loxton--van der Poorten conjecture is proved in \cite{AdBu_Beta} following the approach of \cite{Adamczewski_Bugeaud}. 
Under the same restriction, our method is sufficient, without introducing new ideas, to prove the following result. We recall that a Pisot (resp.\ a Salem) number is a real algebraic integer $>1$, whose complex conjugates lie inside the open unit disc (resp.\ inside the closed unit disc, with at least one of them on 
the unit circle). 

\begin{theo}\label{beta1}
Let $\beta$ be a Pisot or a Salem number. Then, the $\beta$-expansion of a Liouville number cannot be generated by a finite automaton.
\end{theo}

 Moreover, we could provide an explicit upper bound for the irrationality measure of $\xi$. 
Note also that Baker proved in \cite{Baker64} an analog of Theorem \ref{baker}
for number fields. This yields the following result. In the sequel, we will denote by $d_{\beta}(x)$ the $\beta$-expansion of the real number $x$.

\begin{theo}\label{beta2} 
Let $\beta$ be a Pisot or a Salem number, and $\xi$ be a real number that does not lie in $\mathbb Q(\beta)$, and such that 
$$d_{\beta}(\xi)=a_0.a_1a_2\ldots a_n\ldots,$$ 
where ${\bf a}=(a_n)_{n\geq 0}$ is an automatic sequence. Let us assume moreover that the  internal sequence associated with ${\bf a}$ begins in an overlap (see Section \ref{auto} for a definition). Then,   $\xi$ is either a $S$-number or a $T$-number.
\end{theo}

However, new ideas are really needed for dealing with an arbitrary algebraic real number $\beta$. Indeed, in the general case the approximants $\alpha_n$ provided by an automatic $\beta$-expansion (see Section \ref{outline}) are such that the inequality 
$$
\vert \xi-\alpha_n\vert\ll H(\alpha_n)^{-1-\varepsilon}
$$
does not necessarily hold. We recall that the height $H(\alpha)$ of an algebraic number $\alpha$ is defined as the height of its minimal polynomial (see Section \ref{koksma} for a definition). 
Another difficulty also appears if we want to replace the R\'enyi $\beta$-expansion by other representations in base $\beta$ (such as for instance the one arising from the lazy algorithm). Indeed, it is not clear how to obtain the analog of Lemma \ref{dist} for such representations. 

\medskip
 
 We end this Section with a digression on $U$-numbers and some other ge\-neralizations of the approach introduced in this paper.   
As was remarked by LeVeque \cite{Leveque}, the class of $U$-numbers can be divided in infinitely many subclasses of interest; each of them corresponding to the approximation by algebraic numbers whose degree is bounded by a fixed positive integer. The simplest subclass, the class of $U_1$-numbers, exactly corresponds to the Liouville numbers. Given a positive integer $m$, a real number $\xi$ is more generally called a $U_m$-number if 
$$w_m(\xi)=+\infty\;\; \mbox{ and } \;\;w_n(\xi)<+\infty, \mbox{ for }0<n<m.$$ Thus, the set of $U$-numbers is exactly the infinite union on all the positive integers $n$ of the sets formed by the $U_n$-numbers. We refer the reader to Section \ref{koksma} for a definition of $w_n(\xi)$.

Another possible generalization of the present work consists in replacing the $b$-adic expansion with  the continued fraction expansion of real numbers. Then, this naturally leads to consider the notion of automatic continued fraction, that is, to consider real numbers whose continued fraction expansion can be generated by a finite automaton. These real numbers were studied by several authors (see for instance \cite{Adamczewski_BugeaudII} and the references therein).  In particular, it is believed but not yet proved that such numbers are either quadratic or transcendental. This corresponds, in this framework, to  the analog of the Loxton--van der Poorten conjecture. Just as, Becker's conjecture can be translated as follows:  non-quadratic automatic continued fractions are all $S$-numbers. Note that it is obvious that no Liouville number can have an automatic continued fraction expansion, since the latter numbers have by definition bounded partial quotients. Thus, the analog to Theorem \ref{liouville} would be that the set of automatic continued fractions does not contain any $U_2$-numbers. The method introduced in this paper is sufficient to prove a first step toward such a result and would even, in this particular case, provide an explicit upper bound for $w^*_2(\xi)$.

\begin{theo}\label{frac}
Let ${\bf a}=(a_n)_{n\geq 0}$ be an automatic sequence of positive integers. Let us assume that the  first letter of the internal sequence associated with ${\bf a}$ appears at least twice (see Section \ref{auto} for a definition). Then, the automatic continued fraction   
$$\xi:=[a_0,a_1,a_2,\ldots]$$
is not a $U_2$-number.
\end{theo}
 
 As an example, let us consider a uniform morphism $\sigma$ defined from the monoid $\{a,b\}^*$ into itself, where $a$ and $b$ are two distinct positive integers. Let us assume that the sequence ${\bf a}=(a_n)_{n\geq 0}$ is a non-eventually periodic fixed point for $\sigma$. Then, the real number 
 $\xi:=[a_0,a_1,a_2,\ldots]$ 
is transcendental (this is proved in \cite{Adamczewski_BugeaudII}) and is not a $U_2$-number as a consequence of Theorem \ref{frac}. 

Note that it could also be interesting to consider continued fractions associated with regular sequences as introduced by Allouche and Shallit \cite{ALSH}. These sequences provide a natural generalization of automatic sequences for sequences with values lying in an infinite set. 

%%%%%%%%%%%%%%%%%%%%%%%%%%%%%%%%%%%%%%%%%%%%%%
\section{Outlines of the proof of Theorems \ref{beta1}, \ref{beta2}, and \ref{frac}}\label{outline}

We begin with the main steps for proving Theorems \ref{beta1} and \ref{beta2}. 
From now on, $\beta$ denotes a Pisot or a Salem number of degree $l$. 
Let $\xi$ be a real number such that    
$$
d_{\beta}(\xi)=0.a_1a_2\ldots a_n\ldots,
$$
where ${\bf a}=(a_n)_{n\geq 1}$ is an automatic sequence. Note that if the sequence ${\bf a}$ is eventually periodic, then $\xi$ belongs to $\mathbb Q(\beta)$ and cannot be a Liouville number. We keep the notation of the proof of Theorem \ref{mesir}. Then, we first replace the sequence of rational approximants $(p_n/q_n)_{n\geq 0}$ with a  sequence of algebraic numbers $(\alpha_n)_{n\geq 0}$ lying in the same number field $\mathbb Q(\beta)$. We set 
$$
\alpha_n:=0.U_n\overline{V_n},
$$
where $0.U_n\overline{V_n}$ is a possibly inproper expansion in base $\beta$. There thus exists an integer polynomial $P_n(X)$ of degree at most $r_n+s_n$ and such that 
 $$\alpha_n=\frac{P_n(\beta)}{\beta^{r_n}(\beta^{s_n}-1)}\cdot$$ Since $\beta$ is a Pisot or a Salem number of degree $l$, one can use the product formula and the associated notion of height to easily derive that $H(\alpha_n)\ll(r_n+s_n)^{l-1}\beta^{r_n+s_n}$ (see for instance \cite{Adam} where similar arguments are given). Moreover, we have  
\begin{equation*}
\vert \xi-\alpha_n\vert\ll\frac{1}{\beta^{r_n+s_n+t_n}},
\end{equation*}
and thus, for every positive $\varepsilon$, we obtain that
\begin{equation}\label{h}
\vert \xi-\alpha_n\vert\ll \frac{1}{H(\alpha_n)^{1+1/d-\varepsilon}}\cdot
\end{equation}

We now want to apply a classical trick from Diophantine Approximation, as stated in Lemma \ref{53} (see for instance \cite{Leveque, Guting, Bugeaud} for more details).  The proof of this result essentially relies on the following Liouville type inequality (see for instance \cite{Bugeaud}, p. 227). Let $\alpha$ and $\alpha'$ be two distinct algebraic numbers of degree respectively equal to $i$ and $j$. Then, there exists a constant $c$ only depending  on $i$ and $j$ such that 
$$
\vert \alpha-\alpha'\vert >cH(\alpha)^{-j}H(\alpha')^{-i}.
$$  

\medskip

\begin{lem}\label{53}
Let $\xi$  be a real number. Let us assume that there exists an infinite sequence $(\alpha_n)_{n\geq 1}$ of algebraic numbers of degree $r$ such that:

\medskip

\begin{itemize}
\item[\rm (i)] $H(\alpha_n)<H(\alpha_{n+1})<H(\alpha_n)^s$;

\smallskip

\item[\rm (ii)]  $H(\alpha_n)^{-\eta'}<\vert \xi-\alpha_n\vert<H(\alpha_n)^{-\eta}$.
\end{itemize}
\medskip

\noindent Then, $\xi$ is not a $U_t$-number for every integer $t<\eta$. 
 \end{lem}

Here, we can choose $\varepsilon$ small enough in (\ref{h}) 
to ensure that $\eta=1+1/d-\varepsilon>1$, and we want to prove that $\xi$ is not a Liouville number or, equivalently, that $\xi$ is not a $U_1$-number. It thus remains to prove that the sequence $(\alpha_n)_{n\geq 0}$ enjoys the two following properties:
\medskip

\begin{itemize}
\item $\alpha_n$ is not a too good approximation to $\xi$; this corresponds to the left-hand side of $(ii)$;

\smallskip

\item there exists a subsequence of $(\alpha_n)_{n\geq 0}$ satisfying $(i)$.
\end{itemize}

As in the proof of Theorem \ref{mesir}, both properties are ensured thanks to Lemma \ref{kp} and the following annalog of Lemma \ref{dist}. Since this is not a classical result and some difficulties appear, we choose to give below a complete proof of Lemma \ref{dist'}.

\begin{lem}\label{dist'}
Let $\beta>1$ be a Pisot or a Salem number of degree $l$. Let $U$ and $V$ be two finite words defined over the alphabet $\{0,1,\ldots,\lfloor \beta\rfloor\}$ with length respectively equal to $r$ and $s$, and $\alpha$ be the element of $\mathbb Q(\beta)$ defined by
$$\alpha:=\sum_{n=1}^{+\infty}\frac{a_n}{\beta^n},
$$
 where $(a_n)_{n\geq 1}$ is the eventually periodic sequence with preperiod $U$ and period $V$. Let $\xi$ be a real number, $0<\xi<1$,  such that $d_{\beta}(\xi)=0.b_1b_2\ldots b_n\ldots$, where $(b_n)_{n\geq 1}$ is such that there exists
a positive integer $j>r+s$ satisfying:

\begin{itemize}
\item[\rm (i)] $a_n=b_n$, for $1\leq n< j$;

\smallskip

\item[\rm (ii)] $a_j\not=b_j.$
\end{itemize}
Then, 
$$\left\vert\xi- \alpha\right\vert>\frac{1}{(s+1)^{l-1}\beta^{j+s+l-1}}\cdot$$
\end{lem}

\begin{proof}
The key observation is that for every non-negative integer $r$ 
\begin{equation}
\label{bb}
\sum_{n\geq r+1}\frac{b_n}{\beta^n}<\frac{1}{\beta^{r}}\cdot
\end{equation}
Note that such an inequality is in general not satisfied by arbitrary representations in base $\beta$. However, it holds when one considers the $\beta$-expansion (as a by-product of the fact that the $\beta$-expansion arises from the greedy algorithm). 

We first infer from the fact that the coefficients $b_n$ are non-negative and from (\ref{bb}) that 
$$
\sum_{n=j-s+1}^{j}\frac{b_n}{\beta^n}\leq \sum_{n\geq j-s+1}\frac{b_n}{\beta^n}<\frac{1}{\beta^{j-s}}
$$
and thus
\begin{equation}\label{poly}
b_{j-s+1}\beta^{s-1}+b_{j-s+2}\beta^{s-2}+\ldots+b_j<\beta^s.
\end{equation}
We also recall that by definition of the sequence $(a_n)_{n\geq 1}$, we have 
\begin{equation}\label{per}
a_n=a_{n+js},\; \mbox{ if $n>r$ and $j\geq 0$}.
\end{equation}

We have now to distinguish two cases. We set $a_j=i$ and $b_j=m$, and we first assume that $i>m$.
 Then, $i$ is a positive integer and, the coefficients $a_n$ also being non-negative, we infer from (\ref{per}) that  
$$\begin{array}{rl}\alpha=&\displaystyle\sum_{n=1}^{+\infty}\frac{a_n}{\beta^n}\geq\sum_{n=1}^{j-1}\frac{a_n}{\beta^n}+\frac{a_j}{\beta^j}+\frac{a_{j+s}}{\beta^{j+s}}+\frac{a_{j+2s}}{\beta^{j+2s}}\\ \\
=&\displaystyle
\sum_{n=1}^{j-1}\frac{a_n}{\beta^n}+\frac{i}{\beta^j}+\frac{i}{\beta^{j+s}}+\frac{i}{\beta^{j+2s}}>
\sum_{n=1}^{j-1}\frac{a_n}{\beta^n}+\frac{i}{\beta^j}+\frac{i}{\beta^{j+s}} ,
\end{array}$$
whereas thanks to condition $(i)$ and Inequality (\ref{bb}) we have 
$$\xi=\sum_{n=1}^{+\infty}\frac{b_n}{\beta^n}=\sum_{n=1}^{j-1}\frac{a_n}{\beta^n}+\frac{m}{\beta^j}+\sum_{n\geq j+1}\frac{b_n}{\beta^n}<\sum_{n=1}^{j-1}\frac{a_n}{\beta^n}+\frac{m+1}{\beta^j}\leq \sum_{n=1}^{j-1}\frac{a_n}{\beta^n}+\frac{i}{\beta^j} \cdot$$ 
This yields 
 $$\alpha-\xi>\frac{i}{\beta^{j+s}}\geq\frac{1}{\beta^{j+s}}\cdot$$

Now, let us assume that $m>i$. By condition $(i)$, we have $$b_{j-s+1}b_{j-s+2}\ldots b_{j-1}=a_{j-s+1}a_{j-s+2}\ldots a_{j-1}.$$ Moreover, since $a_j=i\leq m-1=b_j-1$, we deduce from (\ref{poly}) that 
\begin{equation}\label{poly2}
P(\beta):=a_{j-s+1}\beta^{s-1}+a_{j-s+2}\beta^{s-2}+\ldots+a_j<\beta^s-1.
\end{equation}
On the one hand, we infer from (\ref{per}) that  
$$\begin{array}{rl}
\alpha=&\displaystyle\sum_{n=1}^{+\infty}\frac{a_n}{\beta^n}=\sum_{n=1}^{j}\frac{a_n}{\beta^n}
+
\sum_{n=1}^{+\infty}\frac{P(\beta)}{\beta^{j+ns}}\\ \\=&\displaystyle\sum_{n=1}^{j}\frac{a_n}{\beta^n}+ \frac{P(\beta)}{\beta^{j+s}} + \sum_{n=2}^{+\infty}\frac{P(\beta)}{\beta^{j+ns}}\\ \\
=&\displaystyle\sum_{n=1}^{j}\frac{a_n}{\beta^n}+ \frac{P(\beta)}{\beta^{j+s}} + \frac{P(\beta)}{\beta^{j+s}(\beta^s-1)},
\end{array}
$$
and we thus derive from (\ref{poly2}) that 
$$
\alpha
< \sum_{n=1}^{j}\frac{a_n}{\beta^n}+\frac{P(\beta)+1}{\beta^{j+s}}\cdot
$$

On the other hand, we have 
$$
\xi=\sum_{n=1}^{+\infty}\frac{b_n}{\beta^n}\geq \sum_{n=1}^{j}\frac{b_n}{\beta^n}=\sum_{n=1}^{j-1}\frac{a_n}{\beta^n}+\frac{m}{\beta^j}\geq \sum_{n=1}^{j-1}\frac{a_n}{\beta^n}+\frac{i+1}{\beta^j}\cdot
$$
This gives 
\begin{equation}\label{76}
\xi- \alpha> \frac{1}{\beta^j}-\frac{P(\beta)+1}{\beta^{j+s}}=\frac{Q(\beta)}{\beta^{j+s}},
\end{equation}
where $Q(\beta):=\beta^s-P(\beta)-1$. Moreover, (\ref{poly2}) implies that 
$Q(\beta)>0$. Note that $\beta$ is a Pisot or a Salem number of degree $l$ and that $Q$ is an integer polynomial of degree $s$ and of height at most $\lfloor \beta\rfloor$.  Let us denote by $\beta=\beta_1,\beta_2,\ldots,\beta_l$ the algebraic conjugates of $\beta$. Then, $\beta$ being an algebraic integer we have that $Q(\beta)\prod_{j=2}^{l}\vert Q(\beta_j)\vert$ is an integer. Since $Q(\beta)>0$, also $\vert Q(\beta_j)\vert>0$ for $2\leq j\leq l$, and we obtain that $Q(\beta)\prod_{j=2}^{l}\vert Q(\beta_j)\vert>0$. Thus,  $Q(\beta)\prod_{j=2}^{l}\vert Q(\beta_j)\vert\geq 1$. Moreover, for $j\geq 2$, $\vert Q(\beta_j)\vert\leq (s+1)\beta$, since $\vert \beta_j\vert \leq 1$. This yields  
$$Q(\beta)\geq((s+1)\beta)^{-l+1}$$  and we thus infer from (\ref{76}) that $$\xi-\alpha>\frac{1}{(s+1)^{l-1}\beta^{j+s+l-1}}\cdot$$ 
 This ends the proof.
\end{proof}

\medskip

We now consider Theorem \ref{beta2}. 
The number field extension of Baker's result can be stated as follows. Let us assume that there exist $\eta>2$ and an infinite sequence $(\beta_n)_{n\geq 0}$ of distinct algebraic numbers lying in the same number field and such that: $\vert \xi-\beta_n\vert\ll H(\beta_n)^{-\eta}$ and $$\limsup_{n\to\infty}\frac{\log(H(\beta_{n+1}))}{\log(H(\beta_n))}<+\infty\cdot$$Then, $\xi$ is not a $U$-number. 
We proceed as for the proof of Theorem \ref{unombres} to construct a sequence $(\alpha_n)_{n\geq 0}$ such that 
$
\vert\xi-\alpha_n\vert\ll H(\alpha_n)^{-2-\varepsilon},
$ 
for a fixed positive $\varepsilon$. Then, we infer from Lemmata \ref{kp} and \ref{dist'} that we can extract a suitable subsequence from $(\alpha_n)_{n\geq 0}$ to apply Baker's result. 

\bigskip

We end this paper with some hints for the proof of Theorem \ref{frac}. The basic idea consists in replacing the sequence of rationals $(p_n/q_n)_{n\geq 0}$  of the proof of Theorem \ref{mesir} by a sequence $(\alpha_n)_{n\geq 0}$ of quadratic approximations, and then to use Lemma \ref{53} with $t=2$. 

Let $\xi:=[0,a_1,a_2,\ldots]$ be an automatic continued fraction and let us denote by $p_n/q_n$ the $n$-th convergent to $\xi$. Thanks to Theorem \ref{cob}, there exists a pair of morphisms $(\varphi,\sigma)$ and a letter $a$ such that ${\bf a}:=(a_n)_{n\geq 1}=\varphi(\sigma^{\infty}(a))$. Since by assumption the letter $a$ occurs at least twice, there exists 
a finite word $U$ such that $aUa$ is a prefix of $\sigma^{\infty}(a)$. This implies that the sequence ${\bf a}$ begins for every integer $n$ with the word $U_nU'_n$, where $U_n=\varphi(\sigma^{n}(aU))$ and $U'_n=\varphi(\sigma^{n}(a))$. 
We set $\vert U_n\vert=s_n$ and $\vert U'_n\vert=t_n$. Then, our quadratic approximants $\alpha_n$ are defined  as the real numbers having  a periodic continued fraction expansion with period $U_n$, that is, $\alpha_n:=[0,\overline{U_n}]$. 
We first observe that $\alpha_n$ is root of the quadratic polynomial
$$
P_n (X) := q_{s_n-1} X^2 + (q_{s_n} - p_{s_n-1}) X - p_{s_n}
$$
and thus $H(\alpha_n)<q_{s_n}$. 
On the other hand, we have 
$$
\vert\xi-\alpha_n\vert<\frac{1}{q_{s_n+t_n}^2}
$$
since the first $s_n+t_n$ partial quotients of $\alpha$ and $\alpha_n$
are the same (see for instance Lemma 2 of \cite{AdBu_Littlewood}). Moreover, the partial quotients of $\xi$ are bounded and we thus infer from classical inequalities for continuants that 
$q_{s_n+t_n}^2\gg q_{s_n}^{\eta}$ for a fixed $\eta>2$ (see for instance Lemmata 3 and 4 of \cite{AdBu_Littlewood}). We thus obtain 
$$
\vert\xi-\alpha\vert\ll H(\alpha_n)^{-\eta}.
$$

Again, to finish the proof we only have to see that $\alpha_n$ is not a too good approximation to $\xi$  (this correspond to the left-hand side of $(ii)$ of lemma \ref{53}) and that there exists a subsequence of $(\alpha_n)_{n\geq 0}$ satisfying condition $(i)$ of Lemma \ref{53}. These conditions are ensured thanks to Lemma~\ref{kp} and Lemma~\ref{dist''} below, the latter playing the role of Lemma~\ref{dist} in the proof of Theorem~\ref{mesir}. 

\begin{lem}\label{dist''} Let $M$ be a positive real number. Let 
$\alpha = [0, a_1, a_2, \ldots]$ and $\xi =
[0,b_1, b_2, \ldots]$ be real numbers whose partial quotients are
at most equal to $M$. 
Assume that there exists a positive integer $n$ such that  
$a_i = b_i$ for any $i=1, \ldots, n$ and $a_{n+1}\not=b_{n+1}$. 
Then, we have
$$
\vert\xi- \alpha\vert \ge \frac{1}{(M+2)^3 q_n^2},
$$ 
where $q_n$ denotes the denominator
of the $n$-th convergent to $\alpha$.
\end{lem}

This latter result corresponds to Lemma 5 in \cite{AdBu_Littlewood}. We thus refer the reader to that paper for a complete proof. 
%%%%%%%%%%%%%%%%%%%%%%%%%%%%%%%%%%%%%%%%%%%%%%%%

%%%%%%%%%%%%%%%%%%%%%%%%%%%%%%%%%%%%%%%%%%%%

\vspace{0.7 cm}

\noindent Boris Adamczewski  

\noindent  CNRS, Institut Camille Jordan

\noindent   Universit\'e Claude Bernard Lyon 1

\noindent   B\^at. Braconnier, 21 avenue Claude Bernard

\noindent   69622 Villeurbanne Cedex   

\noindent (FRANCE) 
\vskip2mm

\noindent {\tt Boris.Adamczewski@math.univ-lyon1.fr}

\bigskip

 \noindent Julien Cassaigne

\noindent CNRS, Institut de Math\'ematiques de Luminy

 \noindent Campus de Luminy, Case 907

 \noindent 13288 Marseille Cedex 9

\noindent (FRANCE)
\vskip 2mm

\noindent {\tt Julien.Cassaigne@iml.univ-mrs.fr}

%%%%%%%%%%%%%%%%%%%%%%%%%%%%%%%%%%%%%%%%%%%%%%%%
\end{document}